# Accurate and Infinite Prime Prediction from Novel Quasi-Prime Analytical Methodology

**Robert E. Grant[1], Talal Ghannam[2]**

[1]Strathspey Crown Holdings, Crown Sterling. Newport Beach, California, USA.
[2]Crown Sterling, Newport Beach, California, USA.

**Abstract**
It is known that prime numbers occupy specific geometrical patterns or moduli when numbers from 1 to ∞ are distributed around polygons having sides that are integer multiple of number 6. In this paper, we will show that not only prime numbers occupy these moduli, but non-prime numbers sharing these same moduli have unique 'prime-ness' properties. When utilizing digital root methodologies, these non-prime numbers provide a novel method to accurately identify prime numbers and prime factors without trial division or probabilistic-based methods. We will also show that the icositetragon (24-sided regular polygon) is a unique polygon pertaining to prime numbers and their ultimate incidence and distribution.

## I. INTRODUCTION

Prime numbers have always been a source of fascination to mathematicians and scientists alike. Their unique mathematical properties, especially their lack of any deviser but themselves, were and still are the subject of countless investigations and many theorems, most important of which is the fundamental theorem of arithmetic[1,2], which states that any integer >1 can be expressed as the unique product of two or more primes. In this sense, prime numbers can be considered the main block upon which all other numbers are built.

Another aspect of prime numbers that has confounded mathematicians is their lack of an apparent pattern. They appear within the infinite string of numbers in such random fashion that devising a functional equation to correctly predict them, infinitely, is believed by many to be an impossible task[3,4].

Mathematicians would also like to know how many prime numbers exist below a certain number. There are several methods to do so, such as the Riemann prime counting function, which is based on Riemann's zeta-function, one of the most important functions in mathematics[5].

During the past few decades, prime numbers have been playing a more practical role than being mere mathematical curiosities. Most importantly, they have been exploited in encryption and decryption methods, which rest mainly on prime numbers being unique factors of numbers and on their perceived randomness[6].

Nevertheless, prime numbers are not entirely random; in fact, they possess several basic orders or patterns that can be readily identified.
As an example, it is well known that the last digit of any prime number can only be 1, 3, 7 or 9. Numbers 2, 4, 6, and 8 are excluded to ensure that the number is odd, hence is not divisible by the number 2. The number 5 is also excluded lest the number is devisable by 5 and hence not a prime.

In more recent years, scientists and mathematicians have discovered more complex orders hidden within the ostensibly random distribution of prime numbers. One such order has been discovered by two mathematicians from Stanford University, Kannan Soundararajan and Robert Lemke Oliver[7].
As we have seen, one known property of prime numbers is that they must end with either 1, 3, 7 or 9. So, theoretically speaking, if prime numbers occur in a purely random fashion, it shouldn't matter what the last digit of the previous prime is to the next one; each one of the four possibilities should have an equal 25% chance of appearing at the end of the next prime number. However, this was not what they found. After looking for the first 400 billion primes, they found that prime numbers tend to avoid having the same last digit for their immediate prime predecessor. As Dr. Oliver put it, they behave as if they "really hate to repeat themselves." Moreover, they found that primes ending with 3 tend to like being followed by one ending with 9, instead of 1 and 7. (This *numerical affinity* can also be observed using digital root analysis[8], where numbers are differentiated into three groups, referred to as 'triplets' that cluster in infinitely repeating patterns. These groups are [3, 6, 9], [1, 4, 7] and [2, 5, 8].)

Another recent and surprising prime-order was discovered by the theoretical chemist and Princeton professor Salvatore Torquato. Prof. Torquato had the idea of looking at what diffraction pattern prime numbers would produce if they were to be modeled as atom-like particles[9]. He observed that not only do prime numbers create a quasicrystal-like interference pattern, the pattern itself is a self-symmetry fractal that has never been observed before, which they labeled *effectively limit-periodic*. Prof. Torquato stated that: "There is much more order in prime numbers than ever previously discovered".

Some orders of prime numbers are more geometrical than numerical; visible only when numbers are looked at in a non-linear fashion. One such geo-numerical distribution is the so-called Ulam Spiral[10,11] discovered by the mathematician Stanislaw Ulam. In this model, when numbers are distributed in a spiraling pattern, prime numbers are aligned along specific horizontal, vertical, and diagonal lines, as shown below.

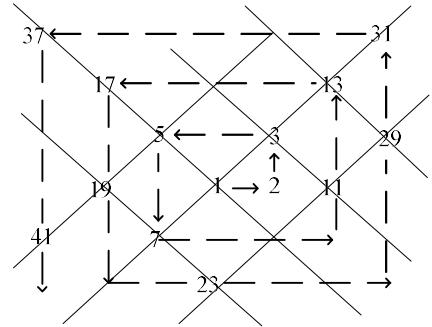

Fig.1: Ulam Spiral where prime numbers are aligned along specific horizontal, vertical and diagonal lines.

Other prime patterns can be found when numbers are distributed around specific concentric polygons. These patterns are at the core of this study as explained in the next section.

## II. POLYGONAL DISTRIBUTION AND QUASI-PRIMES

When numbers from 1 to infinity are distributed within 6 columns, prime numbers will occupy the 1st and 5th columns only (with the exception of 2 and 3), as shown below.

| 1 | 2 | 3 | 4 | 5 | 6 |
|---|---|---|---|---|---|
| 7 | 8 | 9 | 10 | 11 | 12 |
| 13 | 14 | 15 | 16 | 17 | 18 |
| 19 | 20 | 21 | 22 | 23 | 24 |
| 25 | 26 | 27 | 28 | 29 | 30 |
| 31 | 32 | 33 | 34 | 35 | 36 |
| 37 | 38 | 39 | 40 | 41 | 42 |
| 43 | 44 | 45 | 46 | 47 | 48 |

This property emanates from the fact that all prime numbers come in the form $6k \pm 1$, where $k$ is an integer from 0 to infinity.

Consequently, Prime numbers exhibit similar behavior for any distribution of numbers within columns whose number is a multiple of 6. For example, for a given 12-column distribution, prime numbers will fall within the 1st, 5nd, 7th and 11th column, and so on.

Another way to look at these patterns is by distributing numbers around concentric polygons instead of within columns. This distribution has the merit of providing expanded insight regarding the spatial and geometric configuration of prime numbers.

Below are prime distributions around concentric polygons of 24 (4×6) sides, also known as the icositetragon. (For fully illustrated diagram of the 24-wheel, please see appendix C.)

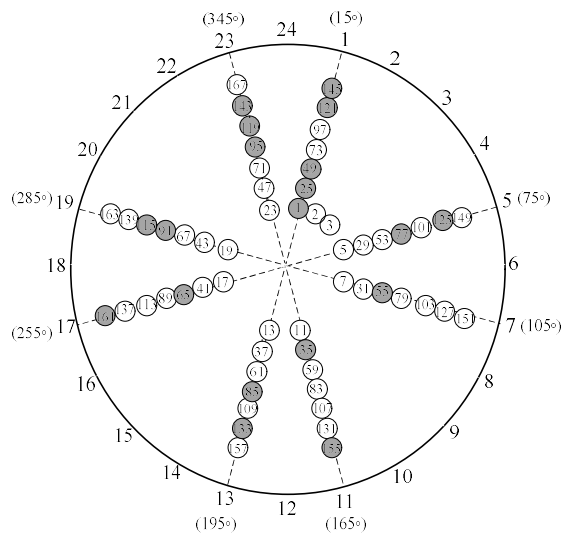

Fig.2: Prime numbers distribution around a 24-sided polygon.

Notice how the moduli prime numbers occupy form shapes that parallel the number of the sides *s* divided by 6. So, for *s* = 24 = 6×4, primes form a double-edged or forked cross occupying moduli 1, 5, 7, 11, 13, 17, 19 and 23[12]. While for *s* = 30 = 6×5, they will form a forked pentagon, and so on. These moduli-based patterns can be easily proven to extend to infinity as follows.

Let us take the 24-fold distribution as an example. Any number on modulus 2 can be written in the form $2 + 24h = 2(1+12h)$, where *h* is any integer from 0 to ∞, and therefore it will always be divisible by 2. Consequently, it cannot be a prime. The same logic applies to Moduli 3, 4, 6, 8, 9, 10, 12, 14, 15, 16, 18, 19, 21, and 22, leaving the rest of the moduli, the *prime moduli*, as the only placeholders for prime numbers.

The above analysis applies to any *s*-sides/columns distribution.

Nevertheless, in the context of prime numbers, the 24-based distribution appears to possess further unique and salient properties. For example, all prime-squared numbers reside in the 1st modulus only. This is not the case for other polygonal distributions.

And when it comes to number 24 specifically, we find that squaring any prime number will always result in a multiple of 24 plus 1 (with the exception of 2 and 3). For example, starting with 5: $5^2 = 25 = 24 + 1$; $19^2 = 24 \times 15 + 1$; $43^2 = 24 \times 77 + 1$, and so on. This is why all prime squares line up along the 1st prime modulus only. (The poof of this can be found in appendix A.)

Ramanujan's Tau function, and its many prime-based conjectures[13], is also based on number 24: $x . \prod_0^\infty (i - x^n)^{24}$.

And as mentioned earlier, it is proven that prime numbers larger than 3 come in the form $6k \pm 1$, in other words, they appear on the sides of numbers that are multiples of 6. Initially, they appear on both sides of multiples of 6, like 5 (6) 7, 11 (12) 13, 17 (18) 19. However, this double pattern breaks up starting at 24, as 25 is not a prime. And again, numbers 2 and 3 are excluded from this pattern. (In fact, these two numbers, 2 and 3, contradict many of the primes properties such that some mathematicians consider them as sub-prime integers. For example, 2 is the only even prime number and 3 is the only prime whose digital root = 3, to be explained later on in the paper. Moreover, notice that all the numbers that are not on the prime-moduli are divisible by 2 and/or 3, while these two numbers will not divide any number on the prime-moduli evenly.)

### III. QUASI-PRIMES, THE Q-GRID AND PREDICTING PRIMES

The above geometrical distributions of prime numbers offer a unique perspective from which we gain important and novel insights. For example, numbers residing along the prime moduli fall within three main categories:

1- Prime numbers,

2- Squared primes, and
3- Quasi-prime numbers.

What we mean by Quasi-primes, are those numbers that are products of primes > or equal to 5 and/or semiprimes. They include numbers such as 55, 175, 245. Notice that they are all odd numbers and are factorized by numbers *residing on the prime-moduli only*. In this sense, they possess certain properties similar to semiprimes; not being actual prime themselves, still, retaining a degree of 'prime-ness'.

However, Quasi-Primes possess the major distinction from semiprimes that they expressly exclude the numbers 2 and 3 as prime factors. In fact, these Quasi-prime numbers are critical in determining infinite primality.

In depth analysis of the 24-distribution reveals that the numbers on the prime moduli are *self-contained within these moduli*. In other words, one can form a multiplication grid, '*Q-grid'*, like the one shown below, where the horizontal and vertical lines consist of prime numbers and semiprimes, e.g. 25, 35, etc. The resultant grid will contain *all the numbers residing along the prime moduli **that are not prime***. Consequently, any number that resides along the prime moduli and is ***not*** contained within the multiplication grid is certainly ***prime***.

|    | 5   | 7   | 11  | 13  | 17  | 19  | 23  | 25  | 29  | 31  | ... |
|----|-----|-----|-----|-----|-----|-----|-----|-----|-----|-----|-----|
| 5  | 25  | 35  | 55  | 65  | 85  | 95  | 115 | 125 | 145 | 155 | ... |
| 7  | 35  | 49  | 77  | 91  | 119 | 133 | 161 | 175 | 203 | 217 | ... |
| 11 | 55  | 77  | 121 | 143 | 187 | 209 | 253 | 275 | 319 | 341 | ... |
| 13 | 65  | 91  | 143 | 169 | 221 | 247 | 299 | 325 | 377 | 403 | ... |
| 17 | 85  | 119 | 187 | 221 | 289 | 323 | 391 | 425 | 493 | 527 | ... |
| 19 | 95  | 133 | 209 | 247 | 323 | 361 | 437 | 475 | 551 | 589 | ... |
| 23 | 115 | 161 | 253 | 299 | 391 | 437 | 529 | 575 | 667 | 713 | ... |
| 25 | 125 | 175 | 275 | 325 | 425 | 475 | 575 | 625 | 725 | 775 | ... |
| 29 | 145 | 203 | 319 | 377 | 493 | 551 | 667 | 725 | 841 | 899 | ... |
| 31 | 155 | 217 | 341 | 403 | 527 | 589 | 713 | 775 | 899 | 961 | ... |
| ... | ... | ... | ... | ... | ... | ... | ... | ... | ... | ... | ... |
|    |     |     |     |     |     |     |     |     |     |     |     |

Notice that the diagonal numbers-line of the grid works as a reflection line with the numbers on top of it being a mirror reflection of those appearing below it. Additionally, the digital root of this line follows a repeated sequence of 174471 174471…

Thus, in theory, the above Q-grid enables us to test the primacy of any number with 100% accuracy and without the need to implement any complicated techniques. It also enables us to predict new primes by simply filtering those numbers residing in the prime-moduli (x (horizontal) and y (vertical) axes) but not appearing as products within the Q-grid.

## IV. PRIME NUMBERS AND DIGITAL ROOT ANALYSIS

All the above patterns of prime numbers, whether numerical or geometrical, can be considered as clearly visible and therefore require minimal analysis to detect.

Other patterns, however, are not readily visible; these require a certain mathematical transformation to be detected, thus requiring the use of digital root analysis.

To calculate the digital root of a specific number, we simply keep summing its individual digits until we are left with a single digit, e.g. $D(137) = D(1+3+7) = D(11) = 2$ (here the letter $D$ denotes taking the digital root of the number).

However, it is when we apply the digital root to sequences of numbers, rather than individual ones, that we are rewarded with some interesting insights.

Take Fibonacci numbers for example; by applying the digital root method to these numbers we find a sequence of digits that repeats itself in a 24-based cycle: [1, 1, 2, 3, 5, 8, 4, 3, 7, 1, 8, 9, 8, 8, 7, 6, 4, 1, 5, 6, 2, 8, 1, 9].

Moreover, when these numbers are distributed around a circle or polygon, every two digits lying on the same diameter will add up to 9 (or a digital root of 9). (This phenomenon further points to the uniqueness of the Icositetragon to prime incidence and distribution).

This is just one example of many fundamental observations we can find in this and many other number sequences[8] such as the Lucas sequence, polygonal numbers, etc.

The digital root is governed by a few mathematical rules, most important of which are the addition and multiplication rules explained as follows.
For $a$, $b$, and $c$ being any rational numbers:
- If $a + b = c$, then
$$D(a + b) = D(a) + D(b) = D(c)$$
- If $a \times b = c$, then
$$D(a \times b) = D(a) \times D(b) = D(c)$$

Of course, the above rules extend to the addition and multiplication of any amount of numbers and is not restricted to two.

Among the simplest prime numbers' patterns one can find using the digital root math is that, and with the exception of number 3, *no prime number can have a digital root of 3, 6, or 9*. Therefore, among the first tests performed on a number to verify its primality should be to calculate its digital root. If the answer is one of these numbers [3, 6, 9], then there is no need to proceed further, even if it is odd and has a last digit of 1, 3, 7, or 9; the number is definitely not a prime.

In fact, any number which has a digital root of 3, 6 or 9 is, in fact, a multiple of 3, and hence cannot be a prime. (The proof of this can be found in appendix B)

This [3, 6, 9] prime fact, along with the general addition and multiplication rules of the digital root analysis, plays an important role in the primary discovery of this research.

## V. PRIMALITY TESTING

There are many methods to test the primality of numbers. The simplest of which is the *trial division* method where if a number $n$ is evenly divisible by any prime number from 2 to $\sqrt{n}$ then it is not prime, and vice versa.

Other methods are considered more probabilistic than exact with the tested numbers subjected to some rigorous criteria, and if they pass, then they are considered most probably prime, such as Fermat Primality test combined with Fibonacci or Lucas primality tests[14, 15].

In this paper, we are proposing an exact and fast method for checking the primality of numbers.
To do so, we start by checking whether a number $n$ passes the initial prime-criteria: being odd, having last-digit of [1, 3, 7 or 9], and digital root not equal to [3, 6, 9]. Once the number passes the above criterion, $n$ is easily tested to whether it falls along any of the prime moduli of the chosen $s$-sides polygon or not. So, for $s = 6$, the prime should fall in either modulus 1 or 5. For $s = 24$, the prime-moduli are: 1, 5, 7, 11, 13, 17, 19, and 23, etc.

By passing all the above criteria, we have reduced the dimensionality of the problem considerably. By requiring the selection of only 8 moduli from the total of 24, we have eliminated 2/3 of the numbers' space that requires search activity. Searching for numbers ending with digits 1, 3, 7, and 9 along these moduli will further reduce the number of possibilities by 1/5. Therefore, the remaining percentage of the numbers' space we need to search in is around 1/3×4/5 = 0.266, less than 30% of the total space.

Once the number passes all the above, we proceed to the final step; to check whether it belongs to the Q-prime grid or not. If the original number $n$ is found within the grid, then it is not a prime. Otherwise, it is a prime by definition.

In case the number $n$ we are testing for primality resides along or close to the reflection line of the Q-grid, then primality confirmation is obvious. We need to look first at $\sqrt{n}$ on both axes of the grid. In case this number is not found we look at those numbers that are immediately bigger and smaller than it. These two numbers will define a small area on the Q-grid where we need to look for this number. Beyond this area, numbers either get larger or smaller than *n*.

On the other hand, if the number *n* is far from the reflection line, then solving the problem requires more complex analysis. Thus, devising indexing algorithm that can search for these numbers without the need to calculate large portions of the Q-grid is an important step in refining and optimizing primality confirmation results.

## VI. PRIME FACTORIZATION

Prime factorization, finding the prime factors of a certain number, is a very difficult task to perform, especially for very large numbers. This is exactly why semi-prime numbers are implemented in many modern encryption techniques[6].

In principle, we can find the prime factors of any semi-prime number, by resorting to the same prime moduli of the Icositetragon, employing both digital root criteria, and Q-grid analysis.

The basic logic in our prime factorization method is similar to that used in primality testing, with some additional restrictions.

For example, having number 1 as the last digit for the number we are trying to factorize informs us regarding the last digits of its two prime factors; namely that they can only be 3 and 7 or both are equal to 1. This will further restrict the size of the Q-grid requiring search analysis.

And because when we search for the prime factors of some semiprime number we need to remove the non-prime numbers from both axes of the Q-grid, like 25, 35, 49 etc., the problem will automatically reduce to simply locating the number in the Q-grid, with its horizontal and vertical projections on both axes being its prime factors.

The advantage of this technique is that it is exact and not based on probability and without requiring division or further mathematical operations.

Similar to the case of primality testing, performing prime factorization using our method requires the implementation of an indexing and search algorithm. The algorithm is programmed based on the search criterion we mentioned earlier along with other supporting conditional statements that exploit the properties and geometry of the Q-grid such that the search is both effective and optimized.

This has already been achieved by the team at Crown Sterling located in Newport Beach, California.

The results were excellent, in terms of accuracy and speed, and were achieved by deploying laptops and desktops with standard specifications and computing power.

## VII. CONCLUSION

By exploiting some of the unique and unfamiliar properties of prime numbers, we were able to explore them in more depth. As a result, a new category of numbers, the Quasi-primes, was identified along with their associated Q-grid.

These discoveries enabled us to tackle the problem of primality and prime factorization with a completely different approach requiring no trial division or probabilistic methodologies, while significantly reducing search scale and volume.

## APPENDIX

A- Proof for any prime number *p*:
$p^2 = k \times 24 + 1$ always.

We know that every prime must come in the form $6k \pm 1$. The factor *k* can be even or odd, so we write it as *2m* for even and *2m+1* for odd. Let us now substitute these in the main equation. We get four forms for *p²* as follows:
$(12m + 1)^2, (12m + 7)^2, (12m - 1)^2, and (12m + 5)^2$

When the squared terms are expanded, we will get four terms that involves 24 times some factor of *m* plus 1, such as $24(6m^2 + m) + 1$ and so on.

B- Proof that no prime number *p* can have a digital root of 3, 6, or 9.

Starting from the proven fact that $p^2 = k \times 24 + 1$, if $p$ has a digital root of 3, 6, or 9 then $D(p) = 9$ always.

On the other side, $D(24) = 6$ and $D(k \times 6) = 3, 6, 9$ always. Now, by adding 1 to the right side, we can never get number 9.

Hence, we can never get 9 on both sides of the equation and consequently no prime number can have a digital root of 3, 6, or 9.

C- The Full 24-Wheel for numbers up to 1008.

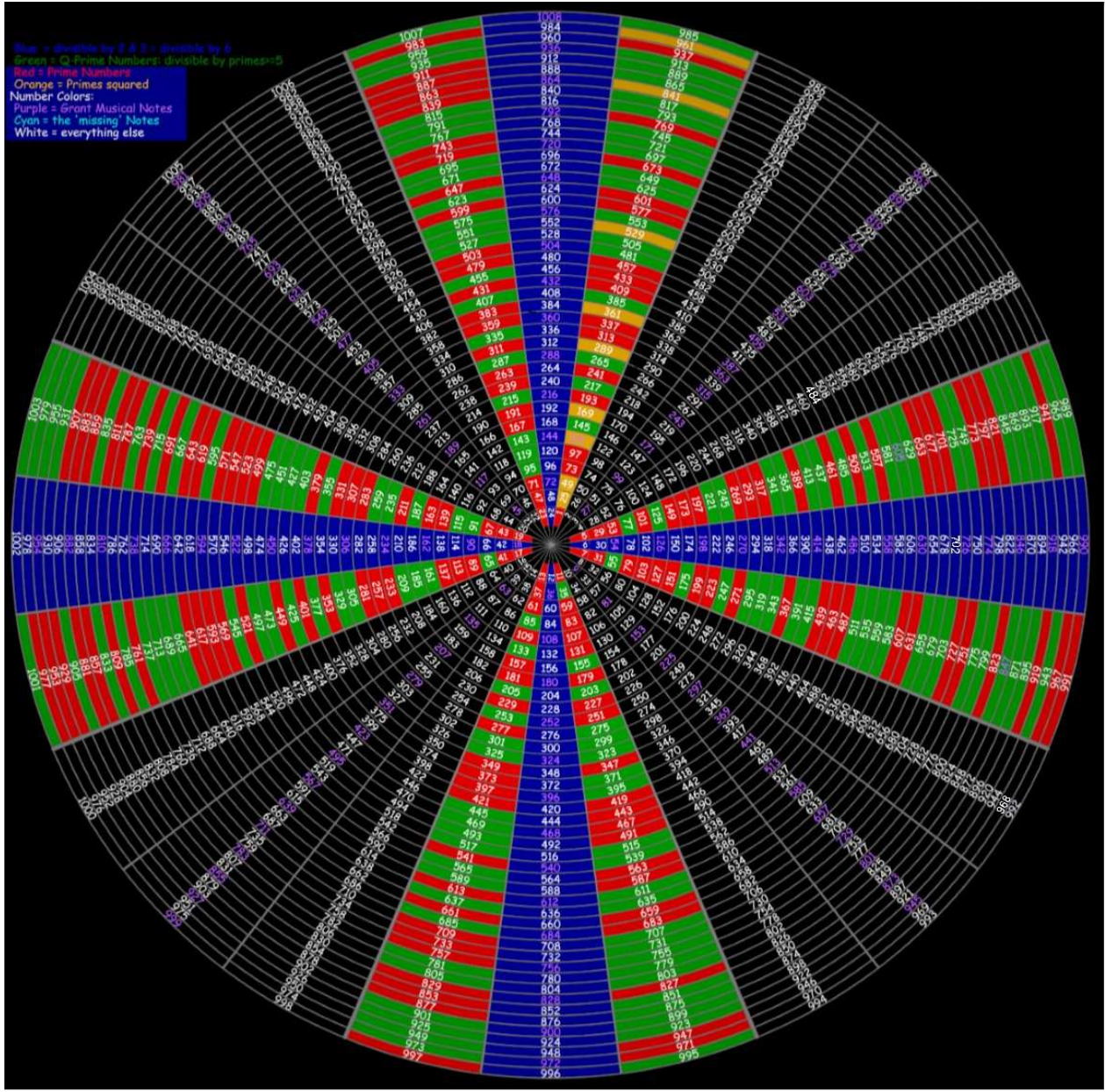